\documentclass{amsart}
\usepackage{amssymb,latexsym}
\theoremstyle{plain}
\newtheorem{theorem}{Theorem}
\newtheorem{corollary}{Corollary}
\newtheorem{proposition}{Proposition}
\newtheorem{lemma}{Lemma}
\theoremstyle{definition}

\newtheorem{remark}{Remark}

\begin{document}

\title[Green and Green-Lazarfeld conjectures for coverings]
{On Green and Green-Lazarfeld conjectures\\ for simple coverings of algebraic curves}
\author{E. Ballico and C. Fontanari}
\address{Department of Mathematics\\
University of Trento\\
Via Sommarive 14\\
38123 Povo (TN), Italy}
\email{ballico@science.unitn.it, fontanar@science.unitn.it}
\thanks{The authors are partially supported by MIUR and GNSAGA of INdAM (Italy).}

\begin{abstract}
Let $X$ be a smooth genus $g$ curve equipped with a simple morphism 
$f: X\to C$, where $C$ is either the projective line or more generally 
any smooth curve whose gonality is computed by finitely many pencils. 
Here we apply a method developed by Aprodu to prove that if $g$ is 
big enough then $X$ satisfies both Green and Green-Lazarsfeld conjectures. 
We also partially address the case in which the gonality of $C$ is computed 
by infinitely many pencils. 
\end{abstract}

\maketitle

\section{Introduction}

Let $X$ be a smooth complex curve of genus $g$. For any spanned 
$L\in \mbox{Pic}(X)$ and all integers $i, j$ let $K_{i,j}(X,L)$
denote the Koszul cohomology groups introduced in \cite{g}. 
Green's conjecture states that $K_{p,1}(X,\omega_X)=0$ if and only 
if $p \ge g-\mathrm{Cliff}(X)-1$, where $\mathrm{Cliff}(X)$ 
is the Clifford index of $X$, while Green-Lazarsfeld conjecture 
(see \cite{GL}, Conjecture (3.7)) predicts that for every line bundle 
$L$ on $X$ of sufficiently large degree $K_{p,1}(X,L)=0$ if and only if 
$p \ge r-\mathrm{gon}(X)+1$, where $r$ is the (projective) dimension of $L$ 
and $\mathrm{gon}(X)$ is the gonality of $X$. 

Both Green and Green-Lazarsfeld conjectures have been verified 
for the general curve of genus $g$ 
(see \cite{V1}, \cite{V2}, \cite{AV}, \cite{A}) and for the general 
$d$-gonal curve of genus $g$ (see \cite{t} for $d \le g/3$, \cite{V1}, 
Corollary 1 on p. 365, for $d \ge g/3$, \cite{AV}, \cite{a0}). 
In particular, \cite{a0} shows that both conjectures are satisfied 
for any smooth $d$-gonal curve verifying a 
suitable linear growth condition on the dimension of Brill-Noether 
varieties of pencils. Such a condition holds for the general $d$-gonal 
curve, but for special curves it turns out to be rather delicate 
(see \cite{m}, Statement T, and \cite{ap}, Proposition 1.3).  

Here we consider the case in which $X$ is a multiple covering. 
Let $h: A\to B$ be a covering of degree $\ge 2$ between
smooth and connected projective curves. The covering $h$ is 
said to be {\it simple} if it does not factor non-trivially, i.e.
for any smooth curve $D$ such that there are morhisms $h_1: A \to D$ 
and $h_2: D \to B$ with $h = h_2\circ h_1$ the morphism $h_2$
is an isomorphism. Every covering of prime order is simple.
By applying \cite{a0}, Theorem 2, and \cite{h}, Theorem 1, 
we are going to prove the following result.

\begin{theorem}\label{i1}
Let $X$ be a smooth genus $g$ curve equipped with a simple morphism 
$f: X\to C$ of degree $m \ge 2$ , where $C$ is a smooth
curve of genus $q$ whose gonality $z$ is computed 
by finitely many $g^1_z$. Assume $g \ge \max \{1+ mq + (m-1)(2mz-5), 
1+ mq + (m-1)(mz-1), 4mz-9, 3mz-6 \}$. Then $X$ satisfies both 
Green and Green-Lazarsfeld conjectures.
\end{theorem}

In the special case $q=0$, the corresponding notion of simple 
linear series is classical (see for instance \cite{ac}) and 
Green's conjecture has already been established for $m \ge 5$ 
in \cite{AF}, Theorem 4.9, by exploiting \cite{c} instead of 
\cite{h}. By the way, for $q=0$ our previous statement simplifies 
as follows. 

\begin{corollary}
Let $X$ be a smooth genus $g$ curve carrying a simple $g^1_m$  
of degree $m \ge 3$. If $g \ge 1+(m-1)(2m-5) \ge 4$ then  
$X$ satisfies both Green and Green-Lazarsfeld conjectures.
\end{corollary}

If instead $q > 0$ and we drop the assumption that the gonality 
of $C$ is computed by finitely many pencils, we obtain with the 
same method the following partial result. 

\begin{proposition}\label{i1bis}
Let $X$ be a smooth genus $g$ curve equipped with a simple morphism 
$f: X\to C$ of degree $m \ge 2$ , where $C$ is a smooth
curve of genus $q \ge 1$ and gonality $z \ge 2$. Assume 
$g \ge 1+ mq + (m-1)(2mz-5)$. If $m=2$ assume also $g \ge 8z-9$. 
Then $K_{p,1}(X,\omega_X)=0$ for any $p \ge g-mz+2$
and $K_{r-mz+2,1}(C, L)=0$ for every line bundle $L$ on $C$ with 
$h^0(C, L)=r+1$ and $\deg(L) \ge 3g$.
\end{proposition}

Finally, if $f: X \to C$ is not simple, then $f= f_s\circ \cdots 
\circ f_1$ with $s\ge 2$ and each $f_i $ a simple covering. 
One could hope to apply Theorem \ref{i1} to each covering $f_i$, 
but the numerical restrictions on the intermediate curves
make such an iterative approach effective only in very few cases.

The authors are grateful to Marian Aprodu and Claire Voisin 
for stimulating e-mail correspondence about Green and 
Green-Lazarsfeld conjectures. 

\section{The proofs}

\begin{remark}\label{i2}
Let $u: X' \to C'$ be a degree $m$ morphism between smooth curves with $X'$ 
of genus $g$ and $C'$ of genus $q$. 
Let $v: X'\to \mathbb {P}^1$ be a degree $x$ morphism such that 
the associated morphism $(u,v): X'\to C'\times \mathbb {P}^1$
is birational onto its image. Then $g \le mq +(m-1)(x-1)$ (Castelnuovo-Severi 
inequality, see for instance \cite{k}, Corollary at p. 26). Notice that $(u,v)$ 
is not birational onto its image if and only if there are a smooth curve $C''$ 
(namely, the normalization of $(u,v)(X')$) and morphisms $w: X' \to C''$, 
$u_1: C'' \to C'$ and $v_1: C'' \to \mathbb {P}^1$ such that $\deg(w) \ge 2$, 
$u = u_1 \circ w$ and $v = v_1 \circ w$. If $u$ is simple, then 
$u_1$ must be an isomorphism and $(u,v)$ is not birational onto 
its image if and only if there is a morphism $\eta = v_1 \circ u_1^{-1}: 
C' \to \mathbb {P}^1$ such that $v = \eta \circ u$. Hence in the set-up 
of Theorem \ref{i1} if $g \ge 1+mq + (m-1)(mz-1)$ then $X$ has gonality $mz$ 
and for every $L\in \mbox{Pic}^{mz}(X)$ such that $h^0(X,L)=2$ there is 
$R\in \mbox{Pic}^z(C)$ such that $h^0(C,R)=2$ and $L \cong f^\ast (R)$.
\end{remark}

\qquad {\emph {Proof of Theorem \ref{i1}.}} By Remark \ref{i2}, 
$X$ has gonality $mz$ and $\dim (W^1_{mz}(X))=0$. M. Aprodu proved 
that $X$ has Clifford index $mz-2$ and satisfies both Green and Green-Lazarsfeld conjectures if 
$\dim (W^1_{mz+t}(X))\le t$ for every integer $t$ such that 
$0 \le t \le g-2mz+2$
(\cite{a0}, Theorem 2). Since the function $x \to \dim (W^1_{x}(X))$ is strictly increasing in the interval $[\mbox{gon}(X),g-1]$,
it is sufficient to prove $\dim (W^1_{g-mz+2}(X))$ $= g-2mz+2$. Assume $\dim (W^1_{g-mz+2}(X)) > g-2mz+2$, i.e. $\dim(W^1_{g-mz+2}(X)) = g -mz-j$ for some integer $j \le mz-3$. 
We have $j \ge 0$ by H. Martens' Theorem (\cite{m}, see for instance \cite{acgh}, IV., Theorem (5.1)). 
Notice also that $g \ge 4j+3$ and $2j+2 \le g-mz+2 \le g-1 -j$, hence a theorem of R. Horiuchi 
yields $\dim (W^1_{2j+2}(X)) = j$ (\cite{h}, Theorem 1). 
Let now $\Gamma$ be any irreducible component of $W^1_{2j+2}(X)$ such
that $\dim (\Gamma )=j$. Since $f$ is simple and
$g-mq > (m-1)(2mz-5) \ge (m-1)(2j+1)$, by Remark \ref{i2} 
there are an integer $y \le \lfloor (2j+2)/m\rfloor$, a non-empty
open subset $\Phi$ of $\Gamma$ and an open subset $\Psi$ of $W^1_y(C)$ 
such that every element of $\Phi$ is the pull-back of an element of 
$\Psi$ plus $2j+2-my$ base points. Thus $j = \dim (\Gamma ) = 
\dim (\Psi )+2j+2-my \le \dim (W^1_y(C))+2j+2-my$.
We have $\Phi \ne \emptyset$, so $\Psi \ne \emptyset$ and $y \ge z$. 
Hence $\dim (W^1_y(C)) \le \dim (W^1_z(C))+2(y-z)$ (\cite{fhl}, Theorem 1). 
Since by assumption $\dim (W^1_z(C))=0$, by putting everything together 
we get $j \le 2(y-z)+2j+2-my \le 2j+2-mz$, i.e. $j \ge mz-2$, 
contradiction.

\qed

The following auxiliary result provides a suitable generalization of 
\cite{a0}, Theorem 2, by repeating almost verbatim the same proof.  

\begin{lemma}\label{main}
Fix an integer $n \ge 1$ and let $C$ be a smooth $d$-gonal curve 
of genus $g$ such that $\dim G^1_{d+m} \le n-1+m$ for all $m$ 
with $n-1 \le m \le g-2d+n+1$. Then $K_{g-d+n,1}(C, \omega_C)=0$
and $K_{r-d+n,1}(C, L)=0$ for every line bundle $L$ on $C$ with 
$h^0(C, L)=r+1$ and $\deg(L) \ge 3g$.      
\end{lemma}

\proof 
Define integers $k, \nu$ as follows: 
\begin{eqnarray}
\label{one} k &=& g-d+n \\
\label{two} \nu &=& 2k-g  
\end{eqnarray}  
and let $X$ be the stable curve obtained from $C$ by identifying 
$\nu+1$ pairs of general points on $C$. In particular, let $p, q$ 
be a pair of points on $C$ identified to a node on $X$. If 
$K_{k,1}(C, \omega_C(p+q))=0$ then according to \cite{AV}, 
Theorem 2.1, for every effective divisor $E$  of degree 
$e \ge 1$ we have $K_{k+e,1}(C, \omega_C(p+q+E))=0$. Thus 
if $L$ is any line bundle on $C$ of degree $x \ge 3g$, 
then $h^0(C,L-\omega_C(p+q)) \ge 1$ and 
$K_{k+x-2g,1}(C, L)=0$. On the other hand, 
by \cite{AV}, Lemma 2.3 and \cite{V1}, p. 367, 
we have $K_{k,1}(C, \omega_C) \subseteq 
K_{k,1}(C, \omega_C(p+q)) \subseteq 
K_{k,1}(X, \omega_X)$, therefore in order 
to prove our statement we may assume $K_{k,1}(X, \omega_X) \ne 0$
and look for a contradiction. By (\ref{two}), $X$ has genus 
$2k+1$, hence by \cite{A}, Proposition 8, there exists a 
torsion-free sheaf $F$ on $X$ with $\deg(F)=k+1$ and $h^0(X,F) \ge 2$. 
Let $s$ with $0 \le s \le \nu+1$ be the number of nodes at which $F$ 
is \emph{not} locally free. If $f: X' \to X$ is the partial 
normalization of $X$ at all such nodes, then $F=f_*(L)$, where 
$L = f^*(F) / \mathrm{Tors}(f^*(F))$ is a line bundle on $X'$ 
with $\deg L = k+1-s$ and $h^0(X',L)=h^0(X,F) \ge 2$. 
By taking the pull-back of $L$ on $C$, we obtain a 
$g^1_{k+1-s}$ not separating $\nu+1-s$ pairs of general points 
on $C$, hence it follows that $\dim G^1_{k+1-s}(C) \ge \nu+1-s$. 

In order to reach a contradiction, assume first 
$0 \le s \le g-2d+2$
(notice that if $g=2r-1$ and $d=r+1$ this case does not occur). 
From (\ref{one}) we obtain   
$k+1-s = d-2d+g+n+1-s$ with $n-1 \le -2d+g+n+1-s \le g-2d+n+1$.
Hence our numerical hypotheses imply that  
$$
\dim G^1_{k+1-s}(C) \le g-2d+2n-s \le \nu-s.
$$

Assume now $s > g-2d+2$.
We claim that also in this case 
$$
\dim G^1_{k+1-s}(C) = \max_r \{ 2(r-1) + \dim W^r_{k+1-s}(C) \}
< \nu+1-s.
$$
Indeed, we have 
\begin{eqnarray*}
\dim W^r_{k+1-s}(C) &\le& \dim W^1_{k+1-s-(r-1)}(C) \le \\
&\le& \dim W^1_d(C) +2(k+1-s-(r-1)-d) \le \\
&\le& 1 +2(k+1-s-(r-1)-d) 
\end{eqnarray*}
where the second inequality is provided by \cite{fhl}, Theorem 1. 
Hence from (\ref{two}) it follows that 
$\dim W^r_{k+1-s}(C) < \nu+1-s-2(r-1)$ for any $r$, as claimed. 

\qed

\qquad {\emph {Proof of Proposition \ref{i1bis}.}}
We argue as in the proof of Theorem \ref{i1} by applying Lemma \ref{main}
with $n=2$ instead of \cite{a0}, Theorem 2. This time we need to prove 
$\dim (G^1_{g-mz+3}(X)) \le g-2mz+4$, so we assume by contradiction 
$\dim (G^1_{g-mz+3}(X))$ $= g-2mz+1-j$ with $j \le mz-4$. 
Once again the numerical hypotheses of \cite{h}, Theorem 1, 
are easily checked, hence we get $j \le \dim (W^1_{z}(C))+2j+2-mz$. 
Since in any case $\dim (W^1_{z}(C)) \le 1$ by \cite{fhl}, Theorem 1, 
we obtain the desired contradiction $j \ge mz-3$.

\qed

\providecommand{\bysame}{\leavevmode\hbox to3em{\hrulefill}\thinspace}

\end{document}